\documentclass[12pt]{amsart}

\numberwithin{equation}{section}

\newtheorem{The}[subsection]{Theorem}
\newtheorem{Lem}[subsection]{Lemma}

\newtheorem{Cor}[subsection]{Corollary}

\theoremstyle{definition}

\theoremstyle{remark}

\newcommand{\field}{{\mathbb{F}}}

\newcommand{\B}{\mathcal {B}}

\newcommand{\V}{{\mathcal {V}}}

\newcommand{\lt}{\mathop{\rm LT}}
\newcommand{\lm}{\mathop{\rm LM}}

\newcommand{\tat}{t\^ete-a-t\^ete}

\title[$SL_2(\mathbb{F}_q)$-invariants]{The Invariants of the Second Symmetric Power Representation of $SL_2(\mathbb{F}_q)$}

\author{Ashley Hobson}
\address{School of Mathematics, Statistics \&  Actuarial Science \\
 \hfil\break\indent University of Kent, Canterbury, CT2 7NF, UK}
\email{ashleyghobson@googlemail.com}

\author{R.\ James Shank}
\address{School of Mathematics, Statistics \&  Actuarial Science \\
 \hfil\break\indent University of Kent, Canterbury, CT2 7NF, UK}
\email{R.J.Shank@kent.ac.uk}

\thanks{The research of the first author was supported by grants from EPSRC}

\date{\today}

\begin{document}
\begin{abstract}
For a prime $p>2$ and $q=p^n$, we compute a finite generating set for the $SL_2(\mathbb{F}_q)$-invariants of the second symmetric power representation,
showing the invariants are a hypersurface and the field of fractions is a purely transcendental extension of the coefficient field. 
As an intermediate result, we show the invariants of
the Sylow $p$-subgroups are also hypersurfaces.

\end{abstract}

\maketitle

\section{Introduction}\label{INTRO}

Consider the generic binary quadratic form over a field $\field$ of characteristic not $2$:
\begin{equation*} a_0X^2+2a_1XY+a_2Y^2.\end{equation*}
Identifying
\begin{equation*} X=\left[\begin{array}{c} 0\\1\end{array}\right] \,\, \mbox{and
} \,\, Y=\left[\begin{array}{c} 1\\0\end{array}\right] \end{equation*}
induces a left action of the general linear group $GL_2(\mathbb{F})$ on the second symmetric power
$$V:={\rm Span}_{\field}[\,Y^2,2XY,X^2\,]$$ and a right action on the dual $V^*={\rm Span}_{\field}[a_2,a_1,a_0].$
For example
$$\begin{array}{cccccc}\sigma_c=\left[ \begin{array}{cc} 1& c\\ 0&1 \end{array}\right] & {\rm acts\ on}\ V^*\ {\rm as}&
\left[ \begin{array}{ccc} 1&2c&c^2\\ 0&1&c\\ 0&0&1\end{array}\right]&&&
\end{array}
$$
with $a_2=[1\ 0\ 0]$, $a_1=[0\ 1\ 0]$, $a_0=[0\ 0\ 1]$.
The action on $V^*$ extends to an action by algebra automorphisms on the symmetric algebra 
$\field[V]=\field[a_2,a_1,a_0]$. For any subgroup $G\leq GL_2(\field)$, we denote the subring of invariant 
polynomials by $\field[V]^G$.

Throughout we assume that $\field$ has characteristic $p>2$, $q=p^n$ and $\field_q\subseteq\field$.
Thus $SL_2(\field_q)\leq GL_2(\field)$. Our primary goal is to describe $\field[V]^{SL_2(\field_q)}$.
Our work generalises and is inspired by L.E. Dickson's solution to the $q=p$ case \cite[Lecture II, \S 8-9]{Dickson}.

Let $P$ denote the subgroup $\{\sigma_c\mid c\in\field_q\}$. $P$ is a Sylow $p$-subgroup of $SL_2(\field_q)$.
The orbit products 
$$\beta:=\prod a_1P=\prod_{c\in\field_q} (a_1+ca_0)=a_1^q-a_0^{q-1}a_1$$
and, for $k\in\field_q$,
$$\gamma_k:=\prod(a_2-ka_0)P=\prod_{c\in\field_q}\left(a_2+2ca_1+(c^2-k)a_0\right)$$
are clearly $P$-invariant.
The discriminant, $\Delta:=a_1^2-a_0a_2$, is a well-known $SL_2(\field_q)$-invariant.
In Section~\ref{PINV}, we show that $\field[V]^P$ is the hypersurface
generated by $a_0,\Delta,\beta,\gamma_0$ subject to the relation
$$\beta^2=a_0^q\gamma_0+\Delta(\Delta^{\frac{q-1}{2}}-a_0^{q-1})^2.$$

Let ${\mathcal Q}$ denote the set of quadratic residues in $\field_q$ and let $\overline{\mathcal Q}$ denote the
set of quadratic nonresidues, i.e., if $\omega$ is a generator for $\field_q^*$, then ${\mathcal Q}$ consists of
the even powers of $\omega$ and $\overline{\mathcal Q}$ consists of the odd powers.
Define 
\begin{eqnarray*}
\Gamma&:=&\prod_{k\in\overline{\mathcal Q}}\gamma_k,\\
B     &:=&\beta\prod_{k\in {\mathcal Q}}\gamma_k,\\
J     &:=&a_0\gamma_0.
\end{eqnarray*}
In Section~\ref{INV}, we show that $\field[V]^{SL_2(\field_q)}$ is the hypersurface generated by
$\Delta, J, \Gamma, B$ subject to a relation of the form
$$B^2=\Delta^q\Gamma^2+J\Phi(\Delta,J,\Gamma)$$
for some polynomial $\Phi$.

Throughout we use the graded reverse lexicographic (grevlex) order with $a_0<a_1<a_2$.
We will see that the given generating sets for $\field[V]^P$ and $\field[V]^{SL_2(\field_q)}$ are
SAGBI bases with respect to this order. A SAGBI basis is the {\bf S}ubalgebra {\bf A}nalogue of
a {\bf G}r\"obner {\bf B}asis for {\bf I}deals. The concept was introduced independently by
Robbiano-Sweedler \cite{rs} and Kapur-Madlener \cite{km}; a useful reference is Chapter~11 of Sturmfels \cite{Sturmfels}
(who uses the term {\it canonical subalgebra basis}). For background material on the invariant theory of finite groups,
see Benson \cite{Benson}, Derksen-Kemper \cite{DK} or Neusel-Smith \cite{NS}.

\section{$P$-invariants}\label{PINV}

For $\omega\in \field_q^*$, the diagonal matrix
$$\begin{array}{cccccc}\rho_{\omega}=\left[ \begin{array}{cc} \omega& 0\\ 0& 1 \end{array}\right] 
& {\rm acts\ on}\ V^*\ {\rm as}&
\left[ \begin{array}{ccc} \omega^2&0&0\\ 0&\omega&0\\ 0&0&1
\end{array}\right].&&&
\end{array}
$$
This motivates the definition of a multiplicative weight function on monomials by
\begin{equation*} \mbox{wt}(a_i)=i. \end{equation*}
Thus for any monomial $\beta$, we have $(\beta)\rho_{\omega}=\omega^{\mbox{wt}(\beta)}\beta$.
Since $\omega^{q-1}=1$, it is convenient to assume that the weight function takes values in
$\mathbb{Z}/(q-1)\mathbb{Z}$.

\begin{Lem} If $f$ is an isobaric polynomial of weight $\lambda$ and $|f P|=|P|$ (i.e., the stabiliser subgroup of $f$ is trivial), 
then $\prod f P$ is isobaric of weight $\lambda$. 
\label{Isobaricprodtransferlem}\end{Lem}

\begin{proof} Note that $P$ is normal in the subgroup of upper-triangular matrices. 
Thus, for $\omega \in \mathbb{F}_q^*$,
\begin{eqnarray*} 
\left(\prod f P\right)\rho_{\omega}
&=&\prod_{\sigma \in P} f\sigma\rho_{\omega} 
=\prod_{\sigma' \in P} f \rho_{\omega}\sigma'\\
&=&\prod_{\sigma' \in P} \omega^{\lambda} f\sigma'
=\omega^{\lambda}\prod f P.\end{eqnarray*} 
Thus $\prod fP$ is isobaric of weight $\lambda$. 
\end{proof}

It is clear that $\Delta$ is isobaric of weight $2$. From the lemma, $\gamma_0$ is isobaric of weight $2$ and $\beta$ is isobaric of weight $1$.
Thus our proposed generators for $\field[V]^P$ are all isobaric.

\begin{Lem} The $P$-invariants $a_0, \Delta, \beta$ and $\gamma_0$ satisfy the relation

\begin{equation*} \beta^2=a_0^q\gamma_0 + \Delta (\Delta^{\frac{q-1}{2}}-a_0^{q-1})^2. \end{equation*}

\label{psylowsyz} \end{Lem}

\begin{proof} Define $\zeta=a_0^q\gamma_0 + \Delta (\Delta^{\frac{q-1}{2}}-a_0^{q-1})^2$.
We first show that $\zeta|_{a_1=0}=0$, which implies $a_1$ divides $\zeta$.

Substituting $a_1=0$ in $\gamma_0$ gives
\begin{eqnarray*} \gamma_0|_{a_1=0} &=& \prod_{t \in \mathbb{F}_q} (t^2a_0+a_2) 
= a_2\prod_{t \in \mathbb{F}_q^*} \left(t^2a_0 + a_2\right) \\
&=& a_2a_0^{q-1}\prod_{t \in \mathbb{F}_q^*} \left(t^2 + \frac{a_2}{a_0}\right) 
= a_2a_0^{q-1}\prod_{s \in{\mathcal Q}}\left(\frac{-a_2}{a_0}-s\right)^2\\
&=&a_2a_0^{q-1} \left(\left(\frac{-a_2}{a_0}\right)^{\frac{q-1}{2}}-1\right)^2=a_2\left(\left(-a_2\right)^{\frac{q-1}{2}}-a_0^{\frac{q-1}{2}}\right)^2.
\label{weirdlabel}\end{eqnarray*}
Thus
$$\zeta|_{a_1=0}
=a_0^q a_2\left(\left(-a_2\right)^{\frac{q-1}{2}}-a_0^{\frac{q-1}{2}}\right)^2+(-a_2a_0)\left((-a_0a_2)^{\frac{q-1}{2}}-a_0^{q-1}\right)^2=0.$$
Therefore $a_1$ divides $\zeta$. However, $\zeta$ is isobaric of weight $2$ and $a_1$ is the only variable of odd weight.
Hence $a_1^2$ divides $\zeta$.

Suppose $a_1$ divides $f\in\field[V]^P$. Then $a_1\sigma_c=a_1+ca_0$ divides $f=f\sigma_c$ for every $c\in \field_q$. Therefore $\beta=\prod a_1P$ divides $f$.
Since $a_1^2$ divides $\zeta$, we see that $\beta^2$ divides $\zeta$. By comparing degrees and lead terms, we conclude that $\beta^2=\zeta$, as required.
\end{proof}

\begin{Lem}\label{Pfield} $\field(V)^P=\field(a_0,\beta,\Delta)$.
\end{Lem}
\begin{proof}
It is easy to verify that $\field[a_0,a_1]^P=\field[a_0,\beta]$
(see, for example, \cite[Theorem 3.7.5]{DK}).
Since $\Delta$ has degree $1$ as a polynomial in $a_2$, applying \cite[Theorem 2.4]{CC} gives
$\field(V)^P=\field(a_0,a_1)^P(\Delta)=\field(a_0,\beta,\Delta)$ (see also  \cite{Kang}).
\end{proof}

\begin{Lem}\label{Phsop} $\{a_0,\Delta,\gamma_0\}$ is a homogeneous system of parameters.
\end{Lem}
\begin{proof} Using grevlex with $a_0<a_1<a_2$, the lead monomials are $a_0$, $a_1^2$ and $a_2^q$.
Thus $(a_0,\Delta,\gamma_0)\field[V]$ is a zero-dimensional ideal and $\{a_0,\Delta,\gamma_0\}$ is a 
homogeneous system of parameters.
\end{proof}

\begin{The} $\B:=\{a_0,\Delta,\beta,\gamma_0\}$ is a generating set, in fact a SAGBI basis, for $\field[V]^P$.
\end{The}

\begin{proof} Let $R$ denote the algebra generated by $\B$. Using grevlex with $a_0<a_1<a_2$, 
there is a single non-trivial t\^{e}te-\`{a}-t\^{e}te,
$\beta^2-\Delta^q$, which, using the relation given in Lemma~\ref{psylowsyz}, subducts to $0$. Thus $\B$ is a SAGBI basis for $R$.

Using Lemmas~\ref{Pfield} and \ref{Phsop}, $\field[V]^P$ is an integral extension of $R$ with the same field of fractions.
Thus to show $R=\field[V]^P$, it is sufficient to show that $R$ is normal, i.e., integrally closed in its field of fractions.
Unique factorisation domains are normal; therefore it is sufficient to show $R$ is a UFD.

Using the relation, we see that $R[a_0^{-1}]=\field[a_0,a_0^{-1}][\Delta,\beta]$, with $a_0,\Delta,\beta$ algebraically independent. 
Thus $R[a_0^{-1}]$ is a UFD.
It follows from \cite[Theorem~20.2]{Matsumura} (or \cite[Lemma~6.3.1]{Benson}) that if $a_0R$ is a prime ideal, $R$ is a UFD.

Suppose $f,g\in R$ with $fg\in a_0R$. Since $R$ is graded, we may assume $f$ and $g$ are homogeneous.
Clearly $a_0\field[V]$ is prime. Therefore, without loss of generality, we may assume $f\in a_0\field[V]$.
Hence the lead monomial $\lm(f)$ is divisible by $a_0$.
$\B$ is a SAGBI basis for $R$ and $f\in R$. Thus $f$ subducts to $0$. Using the grevlex order with $a_0$ small,
every monomial of degree ${\rm deg}(f)$, less than $\lm(f)$, is divisible by $a_0$. Thus at each stage of the subduction,
 there is a factor of $a_0$. Hence $f\in a_0R$ and $a_0R$ is prime.
\end{proof}

\section{$SL_2(\mathbb{F}_q)$-invariants}\label{INV}

The group element
$$\begin{array}{cccccc}\tau=\left[ \begin{array}{rc} 0& 1\\ -1& 0 \end{array}\right] 
& {\rm acts\ on}\ V^*\ {\rm as}&
\left[ \begin{array}{crc} 0&0&1\\ 0&-1&0\\ 1&0&0
\end{array}\right].&&&
\end{array}
$$
It is well-known and easily verified that $\{\tau\}\cup P$ generates $SL_2(\field_q)$.
Thus to show that $f\in\field[V]^P$ is $SL_2(\field_q)$-invariant, it is sufficient to show $(f)\tau=f$.

\begin{Lem}
$J$, $\Gamma$ and $B$ are $SL_2(\field_q)$-invariant.
\end{Lem}

\begin{proof} By construction, $J$, $\Gamma$ and $B$ are $P$-invariant. 
A relatively straightforward calculation shows that each of these polynomials is fixed by $\tau$ 
and is therefore $SL_2(\field_q)$-invariant. It is perhaps more instructive to note that
$SL_2(\field_q)$ permutes the lines in $V^*$ 
and that
each of $J$, $\Gamma$, and $B$ is a projective orbit product. 
For example, the stabiliser of the line $a_0\field_q$ has order $q(q-1)$ and $J$ is a product of $q+1$
linear factors, one taken from each line in the orbit of $a_0\field_q$. Similarly, the
stabiliser of $a_1\field$ has order $2(q-1)$ and $B$ is the product of $q(q+1)/2$ linear 
factors, each representing a line in the orbit. The linear factors of $\Gamma$ are of the form
$a_2+2ca_1+(c^2-k)a_0$ for $c\in\field_q$ and $k\in \overline{\mathcal Q}$. Applying $\tau$ gives
\begin{eqnarray*}
(a_2+2ca_1+(c^2-k)a_0)\tau&=&a_0-2ca_1+(c^2-k)a_2\\
&=&(c^2-k)\left(a_2+2a_1\frac{-c}{c^2-k}+a_0\frac{1}{c^2-k}\right).
\end{eqnarray*}
However
$$\frac{1}{c^2-k}=\left(\frac{-c}{c^2-k}\right)^2-\frac{k}{(c^2-k)^2}$$
with $k/(c^2-k)^2\in \overline{\mathcal Q}$. Thus $\tau$ permutes the lines in $V^*$ corresponding
to the linear factors of $\Gamma$.

Since $SL_2(\field_q)$ acts on $J$, $\Gamma$ and $B$ by permuting the linear factors, up to scalar multiplication,
the action on each of these polynomials is by a multiplicative character. However, for $q\not\in\{2,3\}$,
$SL_2(\field_q)$ is simple (see, for example \cite[4.5]{orange}); hence the character is trivial and the polynomials are invariant.
The case $q=2$ is inconsistent with our hypothesis ${\rm char}(\field)>2$. The $q=3$ case was covered by Dickson's work \cite{Dickson} 
(and can be easily verified by computer). 
\end{proof}

\begin{Lem} $\{\Delta, J, \Gamma \}$ is a homogeneous system of parameters for $\field[V]^{SL_2(\field_q)}$.
\end{Lem}

\begin{proof}
Without loss of generality, we may assume $\field$ is algebraically closed.
We will show that the variety associated to $(\Delta, J, \Gamma)\field[V]$,
say $\V$, consists of the zero vector.

Suppose $v\in\V$. Since $J(v)=0$, there exits $g\in SL_2(\field_q)$ such that
$a_0g(v)=0$. Replacing $v$ with $g(v)$ if necessary, we may assume $a_0(v)=0$.
Thus $\Delta(v)=a_1^2(v)$, giving $a_1(v)=0$. Since $\Gamma\in a_2^{q(q-1)/2}+(a_0,a_1)\field[V]$,
we have $\Gamma(v)=a_2^{q(q-1)/2}(v)$, giving $a_2(v)=0$.
\end{proof}

Define $A:=\field[\Delta, J, \Gamma ]$.

\begin{Cor}\label{rank-cor}
$\field[V]^{SL_2(\field_q)}$ is a free $A$-module of rank $2$.
\end{Cor}
\begin{proof} It is well known that the ring of invariants of a $3$ dimensional representation is Cohen-Macaulay (see \cite[3.4.2]{DK} or \cite[5.6.10]{NS}), i.e., a free module
over any homogeneous system of parameters (hsop). For a faithful action, the rank is given by the order of the group divided by the product of the degrees
of the elements in the hsop (see \cite[3.7.1]{DK} or \cite[5.5.8]{NS}). $SL_2(\field_q)$ acts on $V$ with kernel generated by $-I$ and 
$${\rm deg}(\Delta){\rm deg}(J){\rm deg}(\Gamma)=2(q+1)\frac{q(q-1)}{2}=|SL_2(\field_q)|.$$
Thus $\field[V]^{SL_2(\field_q)}$ has rank $2$ over $A$.
\end{proof}

\begin{The} \label{sl-gen}
$\field[V]^{SL_2(\field_q)}$ is generated by $\Delta$, $J$, $\Gamma$ and $B$ subject to a relation of the form
$$B^2=\Delta^q\Gamma^2+J\Phi(\Delta,J,\Gamma)$$
for some polynomial $\Phi$. Furthermore, this generating set is a SAGBI basis using the grevlex order with $a_0<a_1<a_2$.
\end{The}
\begin{proof}
For any $f\in\field[V]$ we can write $f=f_e+f_o$ were $f_e$ is a sum of terms of even weight and $f_o$ is a sum of terms of odd weight.
If $f\in\field[V]^P=\field[a_0,\Delta,\gamma_0]\oplus \beta\field[a_0,\Delta,\gamma_0]$, then $f_e\in\field[a_0,\Delta,\gamma_0]$
and $f_o\in \beta \field[a_0,\Delta,\gamma_0]$. It is clear that $\tau$ preserves weight-parity. Thus if $f$ is $SL_2(\field_q)$-invariant
$(f_e)\tau=f_e$ and $(f_o)\tau=f_o$, giving $f_e,f_o\in\field[V]^{SL_2(\field_q)}$. Every odd-weight term is divisible by $a_1$. Hence, every
odd-weight $SL_2(\field_q)$-invariant is divisible by $B$. Thus $\field[V]^{SL_2(\field_q)}=E \oplus BE$, were $E$ denotes the subalgebra of 
even-weight $SL_2(\field_q)$-invariants.  Note that $A\subseteq E$.

Using Corollary~\ref{rank-cor}, there exists a homogeneous $SL_2(\field_q)$-invariant, say $\delta$, such that 
$\field[V]^{SL_2(\field_q)}=A \oplus \delta A$. If ${\rm deg}(\delta)<{\rm deg}(B)$, then $\delta\in E$; hence 
$\field[V]^{SL_2(\field_q)}\subseteq E$, giving a contradiction.
If ${\rm deg}(\delta)>{\rm deg}(B)$, then $B\in A$; hence 
$\field[V]^{SL_2(\field_q)}\subseteq A\subseteq E$, again giving a contradiction.
Therefore   ${\rm deg}(\delta)={\rm deg}(B)$. Comparing Hilbert series, i.e., dimensions of homogeneous components,
we see that $A=E$ and $\field[V]^{SL_2(\field_q)}=A\oplus B A$. Therefore
$\field[V]^{SL_2(\field_q)}$ is generated by $\Delta$, $J$, $\Gamma$ and $B$.

Since $B^2$ has even weight, we have $B^2\in A$. Furthermore, $B^2-\Delta^q\Gamma^2\in \field[V]^{SL_2(\field_q)}$ is zero modulo $a_0$.
Thus $J$ divides   $B^2-\Delta^q\Gamma^2$. The quotient is of even weight and hence is an  element of $A$, say $\Phi(\Delta,J,\Gamma)$.
Therefore $B^2=\Delta^q\Gamma^2+J\Phi(\Delta,J,\Gamma)$.

The lead terms of the generators are $\lt(\Delta)=a_1^2$, $\lt(J)=a_0a_2^q$, $\lt(\Gamma)=a_2^{q(q-1)/2}$ and 
$\lt(B)=a_1^qa_2^{q(q-1)/2}$. Thus the only non-trivial \tat\ is given by $B^2-\Delta^q\Gamma^2$. Hence $\{\Delta, J,\Gamma\}$ is a SAGBI
basis for $A$, $\Phi((\Delta,J,\Gamma)$ subducts to zero using $\{\Delta, J,\Gamma\}$ and $B^2-\Delta^q\Gamma^2$ subducts to zero.
Therefore  $\{\Delta, J,\Gamma, B\}$ is a SAGBI basis for $\field[V]^{SL_2(\field_q)}$.
\end{proof}

\begin{Cor}
Define $$m:=\left\lfloor\frac{1}{2}\left(q+1+q(q-1)/2\right)\right\rfloor\  and\ s:=\left\lfloor\frac{1}{2}\left(1+q(q-1)/2\right)\right\rfloor.$$
Then $\field(V)^{SL_2(\field_q)}=\field( B/\Delta^m, J/\Delta^{(q+1)/2},\Gamma/\Delta^s)$, a purely transcendental extension of $\field$.
\end{Cor}
\begin{proof} Let ${\mathcal F}$ denote the field generated by $\{ B/\Delta^m, J/\Delta^{(q+1)/2},\Gamma/\Delta^s\}$.
Clearly ${\mathcal F}\subseteq \field(V)^{SL_2(\field_q)}$.

Suppose $(q-1)/2$ even. Then $m=\frac{1}{2}(q+1+q(q-1)/2)$ and $s=q(q-1)/4$.
Dividing the homogeneous relation from Theorem~\ref{sl-gen} by $\Delta^{2m-1}$ gives
$$\Delta\left(\frac{B}{\Delta^m}\right)^2=\left(\frac{\Gamma}{\Delta^s}\right)^2
+\left(\frac{J}{\Delta^{(q+1)/2}}\right)\Phi(1,J/\Delta^{(q+1)/2},\Gamma/\Delta^s).$$
Thus $\Delta\in{\mathcal F}$. Therefore $J, \Gamma, B\in {\mathcal F}$, giving
${\mathcal F}=\field(V)^{SL_2(\field_q)}$.

Suppose $(q-1)/2$ is odd. Then $m=\frac{1}{2}(q+\frac{q(q-1)}{2})$ and $s=\frac{1}{2}(\frac{q(q-1)}{2}+1)$
Furthermore $\Gamma$ is of odd degree while $J$ and $\Delta$ are of even degree. Thus $\Gamma$ can not appear in $\Phi$.
Dividing the homogeneous relation from Theorem~\ref{sl-gen} by $\Delta^{2m}$ gives
$$\left(\frac{B}{\Delta^m}\right)^2=\Delta\left(\frac{\Gamma}{\Delta^s}\right)^2
+\left(\frac{J}{\Delta^{(q+1)/2}}\right)\Phi(1,J/\Delta^{(q+1)/2}).$$
Thus $\Delta\in{\mathcal F}$. Therefore $J, \Gamma, B\in {\mathcal F}$, giving
${\mathcal F}=\field(V)^{SL_2(\field_q)}$.

\end{proof}

\ifx\undefined\bysame
\newcommand{\bysame}{\leavevmode\hbox to3em{\hrulefill}\,}
\fi

\end{document}